\documentclass[11pt]{eptcs}
\newcommand{\titlerunning}{$k$-convex polyominoes by semi-perimeter.}
\newcommand{\authorrunning}{Andrew R. Conway \& Anthony J Guttmann}

\usepackage{iftex}

\ifpdf
  \usepackage{underscore}         
  \usepackage[T1]{fontenc}        
\else
  \usepackage{breakurl}           
\fi

\usepackage{enumerate}
\usepackage{graphicx}
\usepackage{amssymb,amsfonts,amsmath}
\usepackage{hyperref}
\def\titlerunning{$k$-Convex Polyominoes by Semi-perimeter}
\def\authorrunning{A.R. Conway \& A.J. Guttmann}
\begin{document}

\title{$k$-Convex Polyominoes by Semi-perimeter.}
\author{Andrew R. Conway
\institute{Fairfield, Vic. 3078, Australia}
\email{andrewcombinatorics@greatcactus.org}
\and
Anthony J. Guttmann
\institute{School of Mathematics and Statistics\\
The University of Melbourne\\
Vic.\ 3010, Australia}
\email{tony.guttmann@gmail.com}
}
\date{}                                           
\maketitle
\begin{abstract}
We give the conjectured solution for the generating function of $k$-convex polyominoes, enumerated by semi-perimeter. The solution was obtained from the analysis of enumeration data that we generated. 
\end{abstract}
\section{Introduction}
{\em Convex polyominoes} refers to polyominoes on the square lattice with the property that their perimeter is equal in length to the perimeter of their minimal bounding rectangle. Alternatively, any  line through the polyomino, either horizontal or vertical, will cut exactly two bonds. Loosely speaking, this means that there are no indentations in the perimeter. 

\begin{figure}[h]
    \centering
    \includegraphics[width=0.25\linewidth]{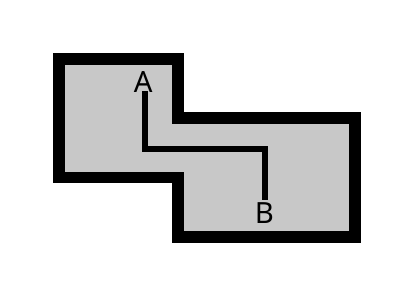}
    \caption{Example 2-convex polyonimo, showing the fewest-corner path between two points A and B.}
    \label{fig:eg}
\end{figure}

A {\em k-convex polyomino} \cite{M22} is a convex polyomino with the additional constraint that any two cells of the polyomino can be joined by a directed {\em internal} NE, NW, SE or SW path with at most $k$ right-angle bends. The path must be parallel to the axes of the polyomino.\footnote{There is an alternative definition of $m$-convex polyominoes in the literature \cite{EGRW}, which refers to a generalisation of convex polyominoes. In that definition, such a polyomino has perimeter equal to the perimeter of its minimal bounding rectangle $+ 2m.$ We are not discussing this class of polyomino here. } See Figure \ref{fig:eg} for an example.

We have developed an algorithm to enumerate $k$-convex polyominoes by both semi-perimeter and area. As the perimeter must be even, it makes sense to refer to the semi-perimeter, otherwise every second coefficient, when enumerating by perimeter, would be zero. Using the perimeter data, we are able to conjecture the exact solution for all $k.$
The description of the enumeration algorithm in appendix \ref{sec:alg} is quite general, and includes enumeration by both perimeter and area. However, in this article we are only subsequently discussing the analysis of the semi-perimeter data.

Our analysis is based on the observation that the two known cases, that of 2-convex \cite{DRS} and fully convex polyominoes \cite{DV} have generating functions of identical structure,
\[
Q(x)=Q_{1}(x)-Q_{2}(x),
\]
where
\[
Q_{1}(x)=\frac{x^{2}p(x)}
{(1-4x)^{2}d(x)},
\qquad
Q_2(x)=\frac{x^{4}q(x)}
{(1-4x)^{3/2}d(x)},
\]
where $p(x),$ $q(x)$ and $d(x)$ are polynomials. We conjecture that this is true for all $k,$ and use our enumerations to identify the polynomials. This cannot be done directly, but requires a non-linear transformation to eliminate the implicit $\sqrt{1-4x}$ in $Q_2(x),$ and transform the result into a rational function. Then a simple Pad\'e approximant suffices to identify that rational function, and, transforming back to the original variable, one can confirm the conjectured form from the known series coefficients that are not implicit in the solution. That is to say, if we have polynomials of degrees that require say, 15 terms for their specification, and we have 25 terms, the extra 10 terms provide resounding  confirmation of our conjectured results.

In this way we conjectured the results for $3 \le k \le 7,$ and from the observed properties of the numerator and denominators, we were able to conjecture simple recurrences, which we were then able to solve.

Our conjectured solution for the semi-perimeter generating function of $k$-convex polyominoes is:
\[
Q_k(x)=Q_{1,k}(x)-Q_{2,k}(x),
\]
where
\[
Q_{1,k}=\frac{-2x^4+x^2(1-6x+7x^2+4x^3) U_k\left (\frac{1-2x}{2x}\right )+x^3(1-4x+8x^2)U_{k-1}\left (\frac{1-2x}{2x}\right )}{(1-4x)^2U_k\left ( \frac{1-2x}{2x} \right )},
\]
and
\[
Q_{2,k}=\frac{4x^{4}T^2_k\left ( \frac{1}{2\sqrt{x}} \right )}{(1-4x)^{3/2}U_k\left ( \frac{1-2x}{2x} \right )}.
\]
Here $T_k(t)$ and $U_k(t)$ are Chebyshev polynomials of the 1st and 2nd kind respectively.

\section{Known results.}

For the full class of convex polyominoes, the generating function is \cite{DV}
\[
C(x)=\frac{x^{2}\left (1-6x+11x^{2}-4x^{3}\right )}{(1-4x)^{2}}
-\frac{4x^{4}}{(1-4x)^{3/2}}.
\]
The number of convex polyominoes of size $n+2$ grows asymptotically as $\frac{n}{8} 4^n$.
For $1$-convex (otherwise known as $L$-convex) polyominoes, the generating function is \cite{CFRR}
 \[
 L(x)=\sum l_n x^n =\frac{x^{2} \left(x^{2}-2 x+1\right)}{2 x^{2}-4 x+1},
 \]
and $ l_n \sim \left ( \frac{1+\sqrt{2}}{4} \right ) \left( 2+\sqrt{2} \right) ^n.$

 For $2$-convex (otherwise known as $Z$-convex polyominoes) the generating function is \cite{DRS}
 \[
    Z(x)= \frac{-x^{2} \left(4 x^{5}-19 x^{4}+48 x^{3}-35 x^{2}+10 x-1\right)}{\left(1-4x\right)^{2} \left(1-3 x\right) \left(1-x\right)} - \frac{x^{4} \left(1-2 x\right)^{2}}{\left(1-4 x\right)^{\frac{3}{2}} \left(1-3 x\right) \left(1-x\right)} \nonumber.
\]
    The number of $Z$-convex polyominoes of size $n+2$ grows asymptotically as $\frac{n}{24} 4^n$.
    
    Prior to this work, the results for $k > 2$ were not known.
    Since the number of $k$-convex polyominoes with $k>2$ is bounded above by the number of convex polyominoes, and bounded below by the number of $(k-1)$-convex polyominoes, it follows that the number of $k$-convex polyominoes ($k>2$) of semi-perimeter $n+2$ grows as $\lambda(k)n 4^n,$ where $\lambda(k)$ is a monotone non-decreasing function growing from $1/24$ to $1/8$ as $k$ goes from 2 to $\infty.$ 
    
    Note that both for $2$-convex and fully convex polyominoes, the structure of the solution is the same:
\[
Q(x)=Q_{1}(x)-Q_{2}(x),
\]
where
\[
Q_{1}(x)=\frac{x^{2}p(x)}
{(1-4x)^{2}d(x)},
\qquad
Q_2(x)=\frac{x^{4}q(x)}
{(1-4x)^{3/2}d(x)}
\]
and $p(x),$ $q(x)$ and $d(x)$ are polynomials. Accordingly,
as a first guess, we suggest that the generating function for $k$-convex polyominoes, enumerated by semi-perimeter is:

\[
Q_k(x)=Q_{1,k}(x)-Q_{2,k}(x),
\]
where
\[
Q_{1,k}(x)=\frac{x^{2}p_k(x)}
{(1-4x)^{2}d_k(x)},
\qquad
Q_2(x)=\frac{x^{4}q_k(x)}
{(1-4x)^{3/2}d_k(x)}.
\]

\section{Conjectured results for $k > 2.$}

To identify this structure from our data is not straightforward. For example, if one asks {\em gfun} to identify the algebraic equation satisfied by the series for 2-convex polyominoes, one requires 30 coefficients to find the equation. As we only have 25 coefficients for $k >2,$ such an approach is clearly not going to be successful.

Instead, we make a change of variable, and write $y=1-\sqrt{1-4x},$ so that $x=\frac{y}{2}-\frac{y^2}{4}.$
Then, for example, the solution for 2-convex polyominoes is
\footnotesize
\[
\frac{y^{10} - 14y^9 + 95y^8 - 376y^7 + 1064y^6 - 2288y^5 + 3632y^4 - 4032y^3 + 2944y^2 - 1280y + 256)(y-2 )^2y^2}{768(y - 1)^4(y^2 - 2y + 4/3)(y^2 - 2y + 4)}
\]
\normalsize
Extracting the factors $\frac{y^2}{(y-1)^4}$ which we know must be present if our assumed form is correct, we are left with a rational function with numerator of degree 12 and denominator of degree 4. This requires 18 series coefficients to identify, simply by constructing a Pad\'e approximant, rather than 30 searching for an algebraic equation from the original series.

\subsection{3-convex polyominoes}
We now apply this approach to our series for 3-convex polyominoes, and readily find the generating function:

\[
Q(x)=Q_1(x)-Q_2(x),
\]
where
\footnotesize
\[
Q_1(x)=\frac{ x^2(8x^6 - 34x^5 + 97x^4 - 110x^3 + 54x^2 - 12x + 1)}
{(1-4x)^{2}(1 - 2x)(2x^2 - 4x + 1)},
\,\,\,
Q_2(x)=\frac{x^{4}\left(1-3x\right)^2}
{(1-4x)^{3/2}(1 - 2x)(2x^2 - 4x + 1)}.
\]
\normalsize
Note that, in the original variable, the numerator of $Q_1(x)$ is of degree 6, the denominator is of degree 3 (we are not worrying about the known terms $\frac{x^2}{(1-4x)^2}$), and $Q_2(x)$ has denominator of degree 2, (and the same numerator). So 14 terms are needed to verify this equation, and we have 25.

Here are the first few terms of the series: \\
\small
\begin{align*}
&(426783233690940x^{24} + 102833217472008x^{23} + 24745444446676x^{22} + 5946372452072x^{21} + 1426784930290x^{20} \\
&+ 341793274268x^{19} + 81735301562x^{18} + 19508859728x^{17} + 4646820820x^{16} + 1104324184x^{15} + 261791420x^{14}\\
& + 61889640x^{13} + 14586482x^{12} + 3426076x^{11} + 801634x^{10} + 186760x^9 + 43302x^8 + 9988x^7 + 2292x^6 + 524x^5 + \\
&120x^4 + 28x^3 + 7x^2 + 2x + 1)x^2.
\end{align*}
\normalsize

The asymptotics are $$a_{n+2} \sim 4^{n}\left (\frac{ n}{16} -\frac{\sqrt{n}}{8\sqrt{\pi}}+\frac{11}{32}+\frac{9}{64\sqrt{n\pi}} \right ).$$ 

Proceeding in this way, we find the following:
\subsection{4-convex by semi-perimeter}
For 4-convex, we find

\[
Q(x)=Q_1(x)-Q_2(x),
\]
where
\[
Q_1(x)=\frac{x^2(-12x^7 + 49x^6 - 178x^5 + 282x^4 - 208x^3 + 77x^2 - 14x + 1)}
{(x^2 - 3x + 1)(5x^2 - 5x + 1)(1 - 4x)^2},
\]
and
\[
Q_2(x)=\frac{x^{4}\left(2x^2 - 4x + 1\right)^2}
{(x^2 - 3x + 1)(5x^2 - 5x + 1)(1 - 4x)^{3/2}}.
\]

Here are the first few terms of the series: \\
\small
\begin{align*}
&(2095409493265028x^{25} + 504284513809564x^{24} + 121184348439516x^{23} + 29075705815146x^{22} +\\
&6964164029736x^{21} + 1664934912802x^{20} + 397228761916x^{19} + 94561501504x^{18} + 22455408404x^{17}\\ 
& + 5318002714x^{16}+ 1255648644x^{15} + 295481944x^{14} + 69273336x^{13} + 16172574x^{12} + 3757972x^{11} +\\
& 868678x^{10} + 199652x^9 + 45608x^8 + 10356x^7 + 2340x^6 + 528x^5 + 120x^4 + 28x^3 + 7x^2 + 2x + 1)x^2.
\end{align*}
\normalsize

The asymptotics are $$a_{n+2} \sim 4^{n}\left (\frac{ 3n}{40} -\frac{\sqrt{n}}{10\sqrt{\pi}}+\frac{21}{80}+\frac{1}{16\sqrt{n\pi}} \right ).$$ 

\subsection{5-convex by semi-perimeter}

For 5-convex, we find

\[
Q(x)=Q_1(x)-Q_2(x),
\]
where
\[
Q_1(x)=\frac{x^2(16x^8 - 84x^7 + 310x^6 - 632x^5 + 644x^4 - 350x^3 + 104x^2 - 16x + 1)}
{(1 - x)(1 - 2x)(1-3x)(1-4x+x^2 )(1 - 4x)^2},
\]
and
\[
Q_2(x)=\frac{x^{4}\left(1-5x+5x^2 \right)^2}
{(1 - x)(1 - 2x)(1-3x)(1-4x+x^2 )(1 - 4x)^{3/2}}.
\]

Here are the first few terms of the series: \\
\small
\begin{align*}
&(9455884509981378x^{26} + 2274678676735172x^{25} + 546358280976036x^{24} + 131014913237644x^{23} + \\
&31360954319330x^{22} + 7492300767612x^{21} + 1786166443572x^{20} + 424837527020x^{19} + 100790612376x^{18} + \\
&23845293592x^{17} + 5623999514x^{16} + 1321924560x^{15} + 309548722x^{14} + 72183568x^{13} + 16755130x^{12} + \\
&3869600x^{11} + 888824x^{10} + 202988x^9 + 46092x^8 +10412x^7 + 2344x^6 + 528x^5 + 120x^4 + 28x^3 + 7x^2 + 2x + 1)x^2.
 \end{align*}
\normalsize
The asymptotics are $$a_{n+2} \sim 4^{n}\left (\frac{ n}{12} -\frac{\sqrt{n}}{12\sqrt{\pi}}+\frac{13}{72}+\frac{5}{288\sqrt{n\pi}} \right ).$$

\subsection{6-convex by semi-perimeter}

For 6-convex, we find

\[
Q(x)=Q_1(x)-Q_2(x),
\]
where
\[
Q_1(x)=\frac{x^2(-20x^9 +127x^8-524x^7+1292x^6-1712x^5+1267x^4-544x^3+135x^2-18x+1)}
{(1-5x+6x^2-x^3)(1-7x+14x^2-7x^3 )(1 - 4x)^2},
\]
and
\[
Q_2(x)=\frac{x^{4}\left(1-4x+x^2 \right)^2(1-2x)^2}
{(1-5x+6x^2-x^3)(1-7x+14x^2-7x^3 )(1 - 4x)^{3/2}}.
\]

Here are the first few terms of the series: \\
\small
\begin{align*}
&(9875153194401768x^{26} + 2371897317966784x^{25} + 568772899988186x^{24} + 136149363839848x^{23} + \\
&32528371648986x^{22} + 7755467989812x^{21} + 1844902796266x^{20} + 437794641248x^{19} + 103609705944x^{18} +\\
& 24448602992x^{17} + 5750559506x^{16} +1347831084x^{15} + 314691982x^{14} + 73165644x^{13} + 16933340x^{12} + \\
&3899788x^{11} + 893464x^{10} + 203604x^9 + 46156x^8 +10416x^7+ 2344x^6 + 528x^5 + 120x^4 + 28x^3 + 7x^2 + 2x + 1)x^2.
\end{align*}
\normalsize

The asymptotics are $$a_{n+2} \sim 4^{n}\left (\frac{ 5n}{56} -\frac{\sqrt{n}}{14\sqrt{\pi}}+\frac{11}{112}-\frac{11}{112\sqrt{n\pi}} \right ).$$ 

\subsection{7-convex by semi-perimeter}

For 7-convex, we find

\[
Q(x)=Q_1(x)-Q_2(x),
\]
where
\footnotesize
\[
Q_1(x)=\frac{x^2(24x^{10} - 198x^9 + 867x^8 - 2476x^7 + 4072x^6 - 3896x^5 + 2251x^4 - 798x^3 + 170x^2 - 20x + 1)}
{(1-8x+20x^2-16x^3+2x^4)(1-4x+2x^2)(1-2x )(1 - 4x)^2},
\]
\normalsize
and
\[
Q_2(x)=\frac{x^{4}\left(7x^3-14x^2+7x-1 \right)^2}
{(1-8x+20x^2-16x^3+2x^4)(1-4x+2x^2)(1-2x )(1 - 4x)^{3/2}}.
\]

Here are the first few terms of the series: \\
\small
\begin{align*}
&(10095807143173876x^{26} + 2421965007900888x^{25} + 580044174368274x^{24} + 138664016164056x^{23} + \\
&33083632385428x^{22} + 7876620943924x^{21} + 1870972607370x^{20} + 443313571324x^{19} + 104755646898x^{18} +\\
& 24681075816x^{17} + 5796404302x^{16} +1356560728x^{15} + 316282360x^{14} + 73439268x^{13} + 16976952x^{12} + \\
&3906036x^{11} + 894228x^{10} + 203676x^9 + 46160x^8 + 10416x^7+ 2344x^6 + 528x^5 + 120x^4 + 28x^3 + 7x^2 + 2x + 1)x^2.
\end{align*}
\normalsize

The asymptotics are $$a_{n+2} \sim 4^{n}\left (\frac{ 3n}{32} -\frac{\sqrt{n}}{16\sqrt{\pi}}+\frac{1}{64}-\frac{23}{128\sqrt{n\pi}} \right ).$$ 

\section{General solution}
\subsection{Observations from the above results}
From the above cases, it appears that the general solution of $k$-convex polyominoes by semi-perimeter to be of the same form as the case for convex polyominoes and 2-convex polyominoes.
\[
Q_k(x)=Q_{1,k}(x)-Q_{2,k}(x),
\]
where
\[
Q_{1,k}(x)=\frac{x^{2}p_k(x)}
{(1-4x)^{2}d_k(x)},
\qquad
Q_2(x)=\frac{x^{4}q_k(x)}
{(1-4x)^{3/2}d_k(x)}.
\]
Here $p_k(x)$ appears to be a polynomial of degree $k+3,$ and $q_k(x)$ similarly appears to be a polynomial of degree $2\lfloor \frac{k}{2}  \rfloor.$ Furthermore, the denominator displays so much regularity that it can be calculated explicitly. It appears to be:
 \[
d_k=1+\sum_{j=1}^k (-1)^j \binom{2k-j+1}{j} x^j.
\]

The polynomials $d_k(x)$ satisfy the linear recurrence\footnote{We thank Simone Rinaldi for this observation}
\[d_k(x) = (1-2x)\,d_{k-1}(x) - x^2\,d_{k-2}(x), \qquad k \ge 2,\]
with initial conditions
$d_0(x)=1, \,\,\, d_1(x)=1-2x.$ And this can be readily solved to yield
$$
d_k(x) 
= \frac{1}{\sqrt{1-4x}} \left(\frac{1-2x+\sqrt{1-4x}}{2} \right)^{k+1}
- \frac{1}{\sqrt{1-4x}} \left( \frac{1-2x-\sqrt{1-4x}}{2} \right)^{k+1}.
$$
We can also write this solution more compactly in terms of Chebyshev polynomials of the second kind, $U_k(t).$
That is to say:
\[
d_k(x)=x^k U_k\left( \frac{1-2x}{2x} \right ).
\]

Next, we note that $q_k(x)$ is a perfect square. That is, $q_k(x)=r_k(x)^2,$ where $r_k(x)$ is obviously of degree $\lfloor \frac{k}{2}  \rfloor.$ 
And we observe from the data that $r_k(x)$ also satisfies a simple recurrence:
\[
r_k(x)=r_{k-1}(x)-x\cdot r_{k-2}(x).\] The initial conditions are $r_0(x)=2$ and $ r_1(x)=1.$
This recurrence can also be solved, and the solution is
\[
r_k(x)=\left ( \frac{1+\sqrt{1-4x}}{2} \right )^k + \left ( \frac{1-\sqrt{1-4x}}{2} \right )^k.
\]
This solution can be expressed in terms of Chebyshev polynomials of the first kind, $T_n(t).$ We find that:
\[
r_k(x)=2x^{k/2}T_k\left (\frac{1}{2\sqrt{x}} \right ).
\]
The numerator of $Q_{1,k}(x),$ that is, $p_k(x),$ caused us some trouble, but we eventually found the recurrence: It is 
\[
p_k(x)=(1-x)p_{k-1}(x)+(x^2-x)p_{k-2}(x)+x^3 p_{k-3}(x),
\]
with initial conditions $p_0(x)=4x^3+5x^2-6x+1,$ $p_1(x)=-16x^3+20x^2-8x+1,$ and $p_2(x)=-4x^5 + 19x^4 - 48x^3 + 35x^2 - 10x + 1.$\\
This too can be solved, with $Q=\sqrt{1-4x},$ to give:
\begin{align*}
p_k(x)=-2x^{k + 2} &+\left (\frac{1 - 2x + Q}{2}\right )^k\left [\frac{1 - 8x + 21x^2 - 18x^3 + 8x^4 + Q(1-6x+7x^2+4x^3)}{2Q}\right ]\\
& +\left (\frac{1 - 2x - Q}{2}\right )^k\left [\frac{-(1 - 8x + 21x^2 - 18x^3 + 8x^4) + Q(1-6x+7x^2+4x^3)}{2Q}\right ].
\end{align*}
This looks a lot like the solution above for $d_k(x),$ and indeed we can write:
\[
p_k(x)=-2x^{k + 2}+x^k\left [(1-6x+7x^2+4x^3)U_k\left ( \frac{1-2x}{2x} \right )+(x-4x^2+8x^3)U_{k-1}\left ( \frac{1-2x}{2x} \right ) \right ]
\]
\subsection{Solution}
Putting all this together, we find the solution for $k$-convex polyominoes by semi-perimeter is conjectured to be:
\[
Q_k(x)=Q_{1,k}(x)-Q_{2,k}(x),
\]
where
\[
Q_{1,k}=\frac{-2x^4+x^2(1-6x+7x^2+4x^3) U_k\left (\frac{1-2x}{2x}\right )+x^3(1-4x+8x^2)U_{k-1}\left (\frac{1-2x}{2x}\right )}{(1-4x)^2U_k\left ( \frac{1-2x}{2x} \right )},
\]
and
\[
Q_{2,k}=\frac{4x^{4}T^2_k\left ( \frac{1}{2\sqrt{x}} \right )}{(1-4x)^{3/2}U_k\left ( \frac{1-2x}{2x} \right )}.
\]
Simone Rinaldi \cite{SR} has shown that, in the limit $k \to \infty,$ the known result for convex polyominoes is recovered. This is further evidence for the correctness of our conjectured result.
\subsection{Asymptotics}
We have calculated the first four terms in the asymptotic expansion of the coefficients for $k <8.$ The general form at leading order and sub-leading order is
\[
a_{n+2,k} \sim 4^{n}\left ( \frac{n(k-1)}{8(k+1)} - \frac{\sqrt{n}}{2\sqrt{\pi}(k+1)} \right ).
\]
Note that in the limit as $k \to \infty$ this gives the correct result, to leading order, for convex polyominoes, viz. $a_{n+2} \sim 4^{n}\left ( \frac{n}{8}\right ),$ though the sub-leading term vanishes in that limit, which it doesn't do in the actual case. That is to say, for convex polyominoes one has $a_{n+2} \sim 4^{n}\left ( \frac{n}{8}- \frac{\sqrt{n}}{2\sqrt{\pi}} \right ).$ 

{\flushleft\textbf{Acknowledgements:} }
We would like to thank Simone Rinaldi and Paolo Massazza for suggesting this problem and for several illuminating discussions. We would also like to thank the University of Melbourne for making available the computing resources enabling these computations to be made.


%
%
%
%
%


\bibliographystyle{eptcs}
\bibliography{Conway}
\appendix

\section{ Description of the algorithm for counting $k$-convex polyominoes.}

\label{sec:alg}

A fairly standard dynamic programming/transfer matrix algorithm is used. 
For explanatory purposes we will start with a description of a well known algorithm
to enumerate column-convex polyominoes. This will then be modified to a slightly more
complex algorithm to enumerate convex polyominoes, and then modified again to
keep track of $k$ for $k$-convexity.

The first two simple algorithms are not necessarily the most efficient known methods of solving
said problems, but are ways that can be extended to the $k$-convex problem in a
relatively straightforward manner.

\subsection{ Simple algorithm 1: Counting column-convex polyominoes.}

A column-convex polyomino on the square lattice has a convexity constraint that only
applies to columns. That is, each column of the polyomino contains a consecutive
range of cells with no holes.
The polyomino can be constructed as a sequence of columns.
To ensure connectivity, each column must have non-zero overlap with the next column.
Suppose one is enumerating by area up to a maximum area $A$.
We will generate a results vector $R$ indexed by the area used, by the following steps

\begin{enumerate}[{\bf1.}]
\item {\bf Initialization.} Start with the leftmost column and consider all columns of length 1 to $A$. Make a list $L$ of
each of these, containing three values:
\begin{itemize}
\item
{\em Signature}: The height of the column (1 to $A$).
\item
{\em Variable Use} : The total area used so far (which will be the same as the height at this point).
\item
{\em Multiplicity} : The number of ways of getting to this point (which will be 1 at this point).
\end{itemize}

\item {\bf Next column.} Next take each element $e \in L$ of the list produced by the previous step and consider what could
go in the next column. 

The next column may be empty - the polyomino is finished. In this case add the multiplicity
of $ e$ to the results vector in the position indicated by the variable use of $e$.

The next column may be non-empty. Consider all possible column lengths and relative positions
of the column relative to the prior column. This will be constrained by the 
maximum area desired and the connectivity constraint. For each of these possibilities,
add a new element to a list $LL$ containing 
\begin{itemize}
\item
{\em Signature}: The height of the new column being added ($1$ to $A$ less prior area use of $e$).
\item
{\em Variable Use}: The total area used so far (which will be the height of this column plus the prior area use of $e$)
\item
{\em Multiplicity}: The number of ways to get to this point (which will be the multiplicity of $e$).
\end{itemize}

If multiple elements end up in $\mathit{LL}$ with the same signature and variable use, then merge them into 
one element of $\mathit{LL}$ with multiplicity equal to the sum of the multiplicities of the merged elements.
If the variable use gets above the maximum area $A$, discard the element.

\begin{figure}[h]
    \centering
    \includegraphics[width=0.3\linewidth]{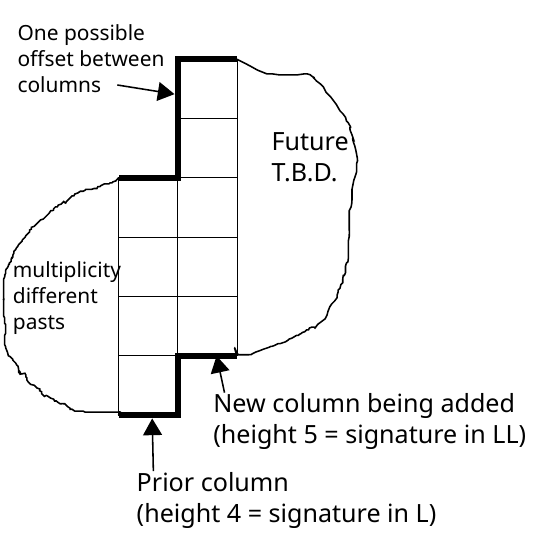}
    \caption{Example addition of a new column. }    \label{fig:IterationStep}
\end{figure}

A picture of an example of one of the possibilities in this step is shown in Figure \ref{fig:IterationStep}.
The past is not relevant to the way columns can connect, 
    other than in the number of different pasts. The future is not determined. The 5 high new column could
    be placed with the top starting anywhere from four units above the top of the prior column to three units below; 
    in this example it is two above. The area variable usage in $L$ is 4 plus the area of the past. The area variable
    usage in $LL$ is 5 plus 4 plus the area of the (now more distant) past. If counting by perimeter, the perimeter 
    usage in $L$ would be 2 plus the perimeter of the past. The perimeter usage in $LL$ would be 5 plus 2 plus
    the perimeter of the (now more distant) past
\footnote{Note: there is a very simple optimization of this algorithm when enumerating just by area 
--- instead of treating the different relative positions 
of the current and prior column separately, just multiply the multiplicity by this number --- that makes 
the algorithm $O(A)$ faster. However it will not be usable for the more complex scenarios described later
so is not relevant for this paper.}.

\item {\bf Completing the algorithm. }
After executing step 2, replace $L$ by $\mathit{LL}$ and repeat step 2. Iterate until the
list of elements is empty. As the variable use will be increased by at least one
for each application of step 2, eventually all elements will exceed the maximum area
constraint and will be discarded resulting in an empty list.

After one iteration of step 2, $R$ will contain all the polyominoes of width 1 (up to the required area), and $LL$ will contain all polyominoes of width exactly 2 (up to the required area), separated by the second column's height.
After two iterations of step 2, $R$ will contain all the polyominoes of width 1 or 2 (up to the required area), and $LL$ will contain all polyominoes of width exactly 3 (up to the required area), separated by the third column's height.
After a sufficient number of iterations, the area (or other) constraint will force $LL$ to be empty.
At this point the results vector contains the number of column-convex polyominoes
for each area up to $A$.

\item  {\bf  Enumerating by other variables.}
The same algorithm can be used to enumerate by other variables - just replace the
variable use computation appropriately. For instance, if counting by perimeter,
the step 1 variable used for a column of height $h$ would be $2h+2$ (including the minimum
closing path), and the addition at step 2 would be 2 horizontal steps and
a number of vertical steps depending upon how the columns overlap.
Similar modifications allow for counting by many other variables, such as the
number of columns, separate vertical and horizontal perimeter, or multiple
simultaneous variables.

\item {\bf Complexity of the algorithm. }
Step 2 is the main complexity.
Suppose there are S possible signatures and V possible variable uses. For enumerating convex poly-
ominoes by area both V and S are O(A).
For step 2, there are O(A) possible next columns and O(A) possible relative positions for each,
leading to O(A2) work for each element of the list There are O(A2) elements of L for a total of O(A4)
work for each iteration of step 2. There are O(A) possible iterations of step 2, for a total of O(A5) time
complexity of the algorithm.
Note that if enumerating by area, as mentioned before the O(A2) work for each element of L can
be optimized to O(A) work by bunching together all the different offsets as they don?t affect anything
needing recording. This reduces the total work from O(A5) to O(A4) in this special case. This does not
work in the more complex scenarios described later, nor when enumerating by perimeter.
Memory use is primarily due to the lists L and LL which are O(A2).

\end{enumerate}

\subsection{ Simple problem 2: Counting convex polyominoes.}

Convex polyominoes are a restriction of column-convex polyominoes
with the additional constraint that each row of the polyomino is a single contiguous range.
This means that adding a new column is more complex than in the previous
algorithm as now columns before the last column matter. However this
can be constrained by adding an extra two Boolean flags to the 
signature. One represents whether a column has ever had a lower top than
the prior column, the other indicates whether a column has ever had a higher bottom
than the prior column. 

In step two, if dealing with an element that has ever had a lower top than
the prior column, then the subsequent column may not have a higher top than the current
column. Similarly, if the element has ever had a higher bottom than the prior
column, then the subsequent column may not have a lower bottom than the current
column.
This increases the number of signatures by a maximum of 4, increasing the
time and memory use by a factor of 4 relative to column-convex polyominoes.

\subsection{ Full problem: $k$-convex polyominoes.}
The same algorithm as for counting convex polyominoes will be used, except
the signature will be extended to include some additional information
that can be used to determine the value of $k$ for each polyomino. When
a polyomino is added to the results vector it will also be indexed by this $k.$

Consider a path $P$ between two points in the polyomino such that $P$ is within
the polyomino, and $P$ has the minimum number of bends to go between the two points.
Due to the convexity constraint of the polyomino, there will never be any reason
for $P$ to double back on itself - that is, go left after going right, or go up
after going down, etc. Assume without loss of generality that $P$ starts from the
leftmost of the two points (if they are in the same column then $P$ is a simple line).
This means we only need to consider paths alternating right and up, and paths alternating
right and down.

The following discussion will focus on tracking the paths going right and up; the
paths going right and down should be treated in an almost identical manner. The maximum
path length (number of corners) will then be the maximum of the two sub-problems.
Each time we add a column, we will keep track of
\begin{itemize}
\item
The maximum number of lines $\mathit{MN}$ in the minimum-corner route  needed to get between any two points so far in the polyomino going up and right.
\item
For each point in the polyomino, the best route to other sites in the polyomino ending in a horizontal line going
  through the vertical border line to the right of the current column. These will be used to update $\mathit{MN}$ when
  the next column is added.
  \end{itemize}
In the following discussion, $B$ is the vertical border line to the right of the current column.

\begin{enumerate}[{\bf 1.}]

\item {\bf   What the best route means.} What information do we need for each point in the polyomino? The best route to a point could potentially
start with either a vertical bond or a horizontal bond. Each vertical bond should then go as far as possible up,
and each horizontal bond should go as far as possible to the right. If you have gone beyond the destination point,
that's fine - going beyond counts as reaching it. So for each point we need
\begin{itemize}
\item
 If we are starting with a vertical bond, then:\\
    $\mathit{VN}$: Number of lines in the path to go through $B$. (This is even)\\
    $\mathit{VH}$: Height at which it goes through $B$.
    \item
If we are starting with a horizontal bond, then:\\
    $\mathit{HN}$: Number of lines in the path to go through $B$. (This is odd)\\
    $\mathit{HH}$: Height at which it goes through $B$.
    \end{itemize}

Note that sometimes it will be impossible to start with a vertical or horizontal bond. In this case, assume such a bond
is of zero length, and the other option will always be used in practice.

This means that each site will be associated with a tuple of 4 numbers ($\mathit{VN},$ $\mathit{VH},$ $\mathit{HN},$ $\mathit{HH}$).

\item {\bf   Eliminating dominated information.} There are some sites whose best route information can be discarded as some other site $S$ dominates them - that is, it is
always harder to get from $S$ to somewhere than from the discardable site. In this case there is no point keeping
track of the discardable site. $S$ is said to dominate the other sites.
In particular, when adding a column, the bottom cell will dominate all higher cells. That is we can ignore all
higher cells. Similarly, if there is a cell to the left of said cell, that will dominate it. The only time
adding a column will necessitate adding a cell to keep track of is when the column starts lower than the prior
column, and even then only the lowest cell needs tracking. More formally,
suppose we have two sites, $\mathit{S1}$ with tuple ($\mathit{VN1},$ $\mathit{VH1},$ $\mathit{HN1},$ $\mathit{HH1}$) and $\mathit{S2}$ with tuple  ($\mathit{VN2},$ $\mathit{VH2},$ $\mathit{HN2},$ $\mathit{HH2}$).
Then $\mathit{S1}$ dominates $\mathit{S2}$ if $\mathit{VN1}\ge \mathit{VN2}$ and $\mathit{VH1}\le \mathit{VH2}$ and $\mathit{HN1}\ge \mathit{HN2}$ and $\mathit{HH1}  \le \mathit{HH2}.$ This potentially allows
removal of other sites, particularly if the polyomino goes through a pinch point.
Note that removing this redundant information does not change the answer at all, but it makes the algorithm
much faster as it vastly reduces the number of different signatures.

\item {\bf   Start state.} For a start column of height $ H,$ the initial state will be a single cell (the bottom one) with
$(\mathit{VN}=2,\,\,\mathit{VH}=H-1\,\,,\mathit{HN}=1,\,\,\mathit{HH}=0).$

Note that the $\mathit{VN}=2$ limit start means the algorithm cannot distinguish between $L$ convex polyominoes (1 bend) and
the trivial case of polyominoes with zero bends (polyominoes that are a horizontal or vertical line).
If these are wanted separately they are trivial to enumerate analytically.

\item {\bf  Update.} Consider a new column of height $H'$ starting at a height of $\delta$ above the prior column of height $H.$
If the new column starts below the old column $(\delta<0),$ a new cell must be added with the same values
as the start state  $(\mathit{VN}=2,\,\,\mathit{VH}=H-1\,\,,\mathit{HN}=1,\,\,\mathit{HH}=0).$

All prior cells' information must be updated. Both the vertical and horizontal pairs are updated
in the same way. Let $(\mathit{xN},\mathit{xH})$ be either $(\mathit{VN},\mathit{VH})$ or $(\mathit{HN},\mathit{HH}),$ then the following process makes the new state $(\mathit{xN'},\mathit{xH'}):$
\begin{itemize}
\item
 If one can continue horizontally, that is $\mathit{xH}\ge \delta,$ then $\mathit{xN}'=\mathit{xN}$ and $\mathit{xH'}=\min(\mathit{xH}-\delta,\mathit{H'}-1)$ as the base has changed, 
  and one can't go above the max height.
  \item
 Otherwise one will need to add a vertical step in the last column 
  (there can't have already been one, as otherwise $xH$ would be $H-1$ and if $H-1<\delta$ 
  then the new column is disconnected from the current column).
  Change $\mathit{xN'}$ to $\mathit{xN}+2$ and change $\mathit{xH'}$ to $\min(\mathit{H}-1-\delta,\mathit{H'}-1).$
  \end{itemize}

One should check that none of the new states are dominated. Doing an imperfect job of this will
increase time and memory usage but will not produce incorrect results. The program used does
not do a perfect job of this as the check is fiddly and error-prone. Now we update $\mathit{NL}$, the number of lines needed to get anywhere in the polyomino so far. This will be the
maximum of the existing $\mathit{NL}$ value and the best route to everywhere in the new column.

The maximum number of lines needed to get to anywhere in the new column is the maximum of the current value, and, 
for each of the prior cells $ C,$ the minimum number of lines needed to get to the top cell of the new line 
(everywhere else is easy). So for each of the cells being tracked, we
compute the best route to get there being the minimum of the best route starting vertically 
or the best route starting horizontally. The best route for a pair $ (\mathit{xN},\mathit{xH})$ that has already been updated as above will be:
\begin{itemize}
\item
If $\mathit{xH'}=\mathit{H'}-1,$ then one can just get there with $\mathit{xN'}.$
\item
Otherwise one will need to add a vertical step $\mathit{xN'}+1.$
\end{itemize}

\item {\bf  Complexity.}
The number of possible cells that are not dominated is hard to characterise. In practice
the number of states (and thus the time and memory use) seem to grow exponentially with
the size of the enumeration.

In practice, to count all $k$-convex polyominoes of semi-perimeter up to 25 took about 60 hours of CPU time on
a single processor, and some 800GB of physical memory. To count all $k$-convex polyominoes of area up to 65 cells
took about 80 hours on a single processor, and used 2TB of physical memory.

The program can be found at  \url{ https://gitlab.com/andrewconwaymaths/k_convex_polyominoes.}
\end{enumerate}

\end{document}